# On the Bahadur slope of the Lilliefors and the Cramér–von Mises tests of normality

**Miguel A. Arcones**[1]

*Binghamton University*

**Abstract:** We find the Bahadur slope of the Lilliefors and Cramér–von Mises tests of normality.

## 1. Introduction

The simplest goodness of fit testing problem is to test whether a random sample $X_1, \ldots, X_n$ is from a particular c.d.f. $F_0$. The testing problem is:

$$H_0 : F = F_0, \quad \text{versus } H_1 : F \not\equiv F_0.$$

A common goodness of fit test is the Kolmogorov–Smirnov test (see Chapter 6 in [8]; and Section 5.1 in [13]). The Kolmogorov–Smirnov test rejects the null hypothesis for large values of the statistic

$$(1.1) \qquad \sup_{t \in \mathbb{R}} |F_n(t) - F_0(t)|,$$

where $F_n(t) = n^{-1} \sum_{j=1}^{n} I(X_j \leq t)$, $t \in \mathbb{R}$, is the empirical c.d.f. Another possible test is the Cramér–von Mises test, which is significative for large values of the statistic:

$$(1.2) \qquad \int_{-\infty}^{\infty} [F_n(t) - F_0(t)]^2 dF_0(t).$$

Anderson and Darling [1] generalize the previous test by adding a weight function and considering:

$$(1.3) \qquad \int_{-\infty}^{\infty} [F_n(t) - F_0(t)]^2 \psi(F_0(t)) dF_0(t),$$

where $\psi$ is a (nonnegative) weight function. The asymptotic distribution of the statistics in (1.1)–(1.3) can be found in [20].

A natural definition of efficiency of tests was given by Bahadur [5, 6]. Let $\{f(\cdot, \theta) : \theta \in \Theta\}$ be a family of p.d.f.'s on a measurable space $(S, \mathcal{S})$ with respect to a measure $\mu$, where $\Theta$ is a Borel subset of $\mathbb{R}^d$. Let $X_1, \ldots, X_n$ be i.i.d.r.v.'s with values in $(S, \mathcal{S})$ and p.d.f. $f(\cdot, \theta)$, for some unknown value of $\theta \in \Theta$. Let $\Theta_0 \subset \Theta$ and let $\Theta_1 := \Theta - \Theta_0$. Consider the hypothesis testing problem $H_0 : \theta \in \Theta_0$ versus $H_1 : \theta \in \Theta_1$. The level (or significance level) of the test is

$$\sup_{\theta \in \Theta_0} \mathbb{P}_\theta \{\text{reject } H_0\}.$$

---

[1]Department of Mathematical Sciences, Binghamton University, Binghamton, NY 13902, USA, e-mail: `arcones@math.binghamton.edu`

The $p$–value of a test is the smallest significance level at which the null hypothesis can be rejected. Suppose that a test rejects $H_0$ if $T_n \geq c$, where $T_n := T_n(X_1, \ldots, X_n)$ is a statistic and $c$ is a constant. Then, the significance level of the test is

$$(1.4) \qquad H_n(c) := \sup_{\theta \in \Theta_0} P_\theta(T_n \geq c),$$

where $P_\theta$ denotes the probability measure for which the data has p.d.f. $f(\cdot, \theta)$. The $p$–value of the test is

$$(1.5) \qquad H_n(T_n).$$

Notice that the $p$–value is a r.v. whose distribution depends on $n$ and on the specified value of the alternative hypothesis. Given two different tests, the one with smallest $p$–value under alternatives is preferred. Since the distribution of a $p$–value is difficult to calculate, Bahadur (1967, 1971) proposed to compare tests using the quantity

$$(1.6) \qquad c(\theta_1) := -2 \liminf_{n \to \infty} n^{-1} \ln H_n(T_n) \ a.s.$$

where the limit is found assuming that $X_1, \ldots, X_n$ are i.i.d.r.v.'s from the p.d.f. $f(\cdot, \theta_1)$, $\theta_1 \in \Theta_1$. The quantity $c(\theta_1)$ is called the Bahadur slope of the test. Given two tests, the one with the biggest Bahadur slope is preferred. For a review on Bahadur asymptotic optimality see and Nikitin [16]. The Bahadur slopes of the tests in (1.1) and (1.2) can be found in Chapter 2 in [16].

For the statistic in (1.1), it is known (see [6]) that if $F_0$ is a continuous c.d.f., then

$$(1.7) \qquad \lim_{n \to \infty} n^{-1} \ln H_n(T_n) = -G(\sup_{t \in \mathbb{R}} |F(t) - F_0(t)|) \ \ a.s.$$

when the data comes from the c.d.f. $F$, $F \not\equiv F_0$, where
$$(1.8)$$
$$G(a) = \inf_{0 \leq t \leq 1-a} \left( (a+t) \ln(t^{-1}(a+t)) + (1-a-t) \ln((1-t)^{-1}(1-a-t)) \right).$$

In this paper, we will consider the Bahadur slopes of some tests of normality, i.e. given a random sample $X_1, \ldots, X_n$ from a c.d.f. $F$ we would like to test:

$$(1.9) \qquad H_0 : F \text{ has a normal distribution, versus } H_1 : F \text{ does not,}$$

We would like to obtain results similar to the one in (1.7) for several tests of normality. Reviews of normality tests are [12, 21] and [15].

Lilliefors [14] proposed the normality test which rejects the null hypothesis for large values of the statistic

$$(1.10) \qquad \sup_{t \in \mathbb{R}} |F_n(\bar{X}_n + s_n t) - \Phi(t)|,$$

where $\Phi$ is the c.d.f. of the standard normal distribution, $\bar{X}_n := n^{-1} \sum_{j=1}^{n} X_j$ and $s_n^2 := (n-1)^{-1} \sum_{j=1}^{n} (X_j - \bar{X}_n)^2$. This test can be used because the distribution of (1.10) is location and scale invariant. We will consider the test of normality which rejects the null hypothesis if

$$(1.11) \qquad \sup_{t \in \mathbb{R}} |F_n(\bar{X}_n + s_n t) - \Phi(t)| \psi(t),$$

where $\psi : \mathbb{R} \to [0, \infty)$ be a bounded function.

We also consider the test of normality which rejects normality if

$$(1.12) \qquad \int_{-\infty}^{\infty} [F_n(\bar{X}_n + s_n t) - \Phi(t)]^2 \psi(\Phi(t)) d\Phi(t),$$

where $\psi : \mathbb{R} \to [0, \infty)$ satisfies $\int_{-\infty}^{\infty} \psi(F_0(t)) dF_0(t) < \infty$. Notice that the statistics in (1.11) and (1.12) are location and scale invariant.

In Section 2, we present bounds in the Bahadur slope for the statistics in (1.11) and (1.12). Our techniques are based on the (LDP) large deviation principle for empirical processes in [2, 3, 4]. We refer to the LDP to [10] and [9]. The proofs are in Section 4.

In Section 3, we present some simulations of the mean of the *p*–value for several test under different alternatives. The simulations show that Lilliefors test has a high *p*-value. However, the *p*-value of the Anderson–Darling is competitive with other test of normality such as the Shapiro–Wilk test ([19]) and the BHEP test ([11] and [7]).

## 2. Main results

In this section we review some results on the LDP for empirical processes. We determine the rate function of the LDP of empirical processes using Orlicz spaces theory. A reference in Orlicz spaces is [18]. A function $\Upsilon : \mathbb{R} \to \bar{\mathbb{R}}$ is said to be a Young function if it is convex, $\Upsilon(0) = 0$; $\Upsilon(x) = \Upsilon(-x)$ for each $x > 0$; and $\lim_{x \to \infty} \Upsilon(x) = \infty$. Let $X$ be a r.v. with values in a measurable space $(S, \mathcal{S})$. The Orlicz space $\mathcal{L}^{\Upsilon}(S, \mathcal{S})$ (abbreviated to $\mathcal{L}^{\Upsilon}$) associated with the Young function $\Upsilon$ is the class of measurable functions $f : (S, \mathcal{S}) \to \mathbb{R}$ such that $E[\Upsilon(\lambda f(X))] < \infty$ for some $\lambda > 0$. The Minkowski (or gauge) norm of the Orlicz space $\mathcal{L}^{\Upsilon}(S, \mathcal{S})$ is defined as

$$N_{\Upsilon}(f) = \inf\{\lambda > 0 : E[\Upsilon(f(X)/\lambda)] \leq 1\}.$$

It is well known that the vector space $\mathcal{L}^{\Upsilon}$ with the norm $N_{\Upsilon}$ is a Banach space. Define

$$\mathcal{L}^{\Upsilon_0} := \{f : S \to \mathbb{R} : E[\Upsilon_0(\lambda|f(X)|)] < \infty \text{ for some } \lambda > 0\},$$

where $\Upsilon_0(x) = e^{|x|} - |x| - 1$. Let $(\mathcal{L}^{\Upsilon_0})^*$ be the dual of $(\mathcal{L}^{\Upsilon_0}, N_{\Upsilon_0})$. The function $f \in \mathcal{L}^{\Upsilon_0} \mapsto \ln\left(E[e^{f(X)}]\right) \in \mathbb{R}$ is a convex lower semicontinuous function. The Fenchel–Legendre conjugate of the previous function is:

$$(2.1) \qquad J(l) := \sup_{f \in \mathcal{L}^{\Upsilon_0}} \left(l(f) - \ln\left(E[e^{f(X)}]\right)\right), \ l \in (\mathcal{L}^{\Upsilon_0})^*.$$

$J$ is a function with values in $[0, \infty]$. Since $J$ is a Fenchel–Legendre conjugate, it is a nonnegative convex lower semicontinuous function. If $J(l) < \infty$, then:

(i) $l(\mathbf{1}) = 1$, where $\mathbf{1}$ denotes the function constantly 1.

(ii) $l$ is a nonnegative definite functional: if $f(X) \geq 0$ a.s., then $l(f) \geq 0$.

Since the double Fenchel–Legendre transform of a convex lower semicontinuous function coincides with the original function (see e.g. Lemma 4.5.8 in [9]), we have that

$$(2.2) \qquad \sup_{l \in \mathcal{L}^{\Upsilon_0}} (l(f) - J(l)) = \ln E[e^{f(X)}].$$

The previous function $J$ can be used to determine the rate function in the large deviations of statistics. Let $\{X_j\}_{j=1}^{\infty}$ be a sequence of i.i.d.r.v.'s with the distribution of $X$. If $f_1, \ldots, f_m \in \mathcal{L}^{\Upsilon_0}$, then

$$\{(n^{-1}\sum_{j=1}^{n} f_1(X_j), \ldots, n^{-1}\sum_{j=1}^{n} f_m(X_j))\}$$

satisfies the LDP in $\mathbb{R}^m$ with speed $n$ and rate function

$$I(u_1, \ldots, u_m) := \sup_{\lambda_1, \ldots, \lambda_m \in \mathbb{R}} \left( \sum_{j=1}^{m} \lambda_j u_j - \ln E[\exp(\sum_{j=1}^{m} \lambda_j f_j(X))] \right)$$

(see for example Corollary 6.1.16 in [9]). This rate function can be written as

$$\inf \left\{ J(l) : l \in (\mathcal{L}^{\Upsilon_0})^*, l(f_j) = u_j \text{ for each } 1 \leq j \leq m \right\},$$

(see Lemma 2.2 in [4]).

To deal with empirical processes, we will use the following theorem:

**Theorem 2.1** (Theorem 2.8 in [3]). *Suppose that* $\sup_{t \in T} |f(X, t)| < \infty$ *a.s. Then, the following sets of conditions* ((a) *and* (b)) *are equivalent:*

(a.1) $(T, d)$ *is totally bounded, where* $d(s, t) = E[|f(X, s) - f(X, t)|]$.

(a.2) *There exists a* $\lambda > 0$ *such that*

$$E[\exp(\lambda F(X))] < \infty,$$

*where* $F(x) = \sup_{t \in T} |f(x, t)|$.

(a.3) *For each* $\lambda > 0$, *there exists a* $\eta > 0$ *such that* $E[\exp(\lambda F^{(\eta)}(X))] < \infty$, *where* $F^{(\eta)}(x) = \sup_{d(s,t) \leq \eta} |f(x, s) - f(x, t)|$.

(a.4) $\sup_{t \in T} |n^{-1}\sum_{j=1}^{n}(f(X_j, t) - E[f(X_j, t)])| \xrightarrow{\Pr} 0$.

(b) $\{n^{-1}\sum_{j=1}^{n} f(X_j, t) : t \in T\}$ *satisfies the large deviation principle in* $l_{\infty}(T)$ *with speed* $n$ *and a good rate.*

*Besides, the rate function is*

$$I(z) = \inf\{J(l) : l \in (\mathcal{L}^{\Upsilon_0})^*, l(f(\cdot, t)) = z(t), \text{ for each } t \in T\}, z \in l_{\infty}(T).$$

We will consider large deviations when the r.v.'s have a standard normal distribution. We denote $(\mathcal{L}_{\Phi}^{\Upsilon_0}, N_{\Upsilon_0})$ to the Orlicz space, when the distribution of $X$ is a standard normal one. Similarly,

$$(2.3) \qquad J_{\Phi}(l) := \sup_{f \in \mathcal{L}_{\Phi}^{\Upsilon_0}} \left( l(f) - \ln \left( E_{\Phi}[e^{f(X)}] \right) \right), \ l \in (\mathcal{L}^{\Upsilon_0})^*.$$

First, we consider the Bahadur efficiency of the test in (1.11). Next lemma considers the large deviations of the test statistic in (1.11) under the null hypothesis.

**Lemma 2.1.** *Let $\psi : \mathbb{R} \to \mathbb{R}$ be a bounded function. Then, for each $u \geq 0$,*

$$-\inf\left\{ J_\Phi(l) : l \in (\mathcal{L}_\Phi^{\Upsilon_0})^*, \sup_{t \in \mathbb{R}} |x(a + (b - a^2)^{1/2}t) - \Phi(t)|\psi(t) > u, \right.$$

$$l(f_t) = x(t), t \in \mathbb{R}, l(g) = a, l(g^2) = b,$$

$$\left. \text{where } f_t(s) = I(s \leq t), g(s) = s, s \in \mathbb{R} \right\}$$

$$(2.4) \qquad \leq \liminf_{n \to \infty} n^{-1} \ln P_\Phi \left\{ \sup_{t \in \mathbb{R}} |F_n(\bar{X}_n + s_n t) - \Phi(t)|\psi(t) > u \right\}$$

$$\leq \limsup_{n \to \infty} n^{-1} \ln P_\Phi \left\{ \sup_{t \in \mathbb{R}} |F_n(\bar{X}_n + s_n t) - \Phi(t)|\psi(t) \geq u \right\}$$

$$\leq -\inf\left\{ J_\Phi(l) : l \in (\mathcal{L}_\Phi^{\Upsilon_0})^*, \sup_{t \in \mathbb{R}} |x(a + (b - a^2)^{1/2}t) - \Phi(t)|\psi(t) \geq u, \right.$$

$$l(f_t) = x(t), t \in \mathbb{R}, l(g) = a, l(g^2) = b,$$

$$\left. \text{where } f_t(s) = I(s \leq t), g(s) = s, s \in \mathbb{R} \right\}$$

**Theorem 2.2.** *Let $\psi : \mathbb{R} \to \mathbb{R}$ be a bounded function, let*

$$(2.5) \qquad H_n^{\mathrm{Li}}(u) := P_\Phi \left\{ \sup_{t \in \mathbb{R}} |F_n(\bar{X}_n + s_n t) - \Phi(t)|\psi(t) \geq u \right\}, u \geq 0,$$

*and let*

$$G^{\mathrm{Li}}(u) := \inf\left\{ J_\Phi(l) : l \in (\mathcal{L}_\Phi^{\Upsilon_0})^*, \sup_{t \in \mathbb{R}} |x(a + (b - a^2)^{1/2}t) - \Phi(t)|\psi(t) \geq u, \right.$$

$$l(f_t) = x(t), t \in \mathbb{R}, \ l(g) = a, \ l(g^2) = b,$$

$$\left. \text{where } f_t(s) = I(s \leq t), g(s) = s, s \in \mathbb{R} \right\}$$

*Let $\{X_j\}_{j=1}^\infty$ be a sequence of i.i.d.r.v.'s with c.d.f. $F$. Then,*

$$-\lim_{\delta \to 0+} G^{\mathrm{Li}} \left( \sup_{t \in \mathbb{R}} |F(\mu_F + \sigma_F t) - \Phi(t)|\psi(t) + \delta \right)$$

$$(2.6) \qquad \leq \liminf_{n \to \infty} n^{-1} \ln H_n^{\mathrm{Li}} \left( \sup_{t \in \mathbb{R}} |F_n(\bar{X}_n + s_n t) - \Phi(t)|\psi(t) \right)$$

$$\leq \limsup_{n \to \infty} n^{-1} \ln H_n^{\mathrm{Li}} \left( \sup_{t \in \mathbb{R}} |F_n(\bar{X}_n + s_n t) - \Phi(t)|\psi(t) \right)$$

$$\leq -\lim_{\delta \to 0+} G^{\mathrm{Li}} \left( \sup_{t \in \mathbb{R}} |F(\mu_F + \sigma_F t) - \Phi(t)|\psi(t) - \delta \right) \text{ a.s.}$$

*where $\mu_F = E_F[X]$ and $\sigma_F^2 = \mathrm{Var}_F(X)$.*

For the statistic in (1.12), we have similar results:

**Lemma 2.2.** *Let $\psi : \mathbb{R} \to [0, \infty)$ be a function such that $\int_{-\infty}^\infty \psi(F_0(t)) dF_0(t) < \infty$.*

*Then, for each $u \geq 0$,*

$$-\inf\Big\{J_\Phi(l) : l \in (\mathcal{L}_\Phi^{\Upsilon_0})^*,$$

$$\int_{-\infty}^{\infty}[x(a + (b-a^2)^{1/2}t) - \Phi(t)]^2\psi(F_0(t))dF_0(t) > u,$$

$$l(f_t) = x(t), t \in \mathbb{R}, \ l(g) = a, \ l(g^2) = b,$$

$$\text{where } f_t(s) = I(s \leq t), g(s) = s, s \in \mathbb{R}\Big\}$$

(2.7)

$$\leq \liminf_{n\to\infty} n^{-1}\ln P_\Phi\left\{\int_{-\infty}^{\infty}[F_n(\bar{X}_n + s_n t) - \Phi(t)]^2\psi(F_0(t))dF_0(t) > u\right\}$$

$$\leq \limsup_{n\to\infty} n^{-1}\ln P_\Phi\left\{\int_{-\infty}^{\infty}[F_n(\bar{X}_n + s_n t) - \Phi(t)]^2\psi(F_0(t))dF_0(t) \geq u\right\}$$

$$\leq -\inf\Big\{J_\Phi(l) : l \in (\mathcal{L}_\Phi^{\Upsilon_0})^*,$$

$$\int_{-\infty}^{\infty}[x(a + (b-a^2)^{1/2}t - \Phi(t))]^2\psi(F_0(t))dF_0(t) > u,$$

$$l(f_t) = x(t), t \in \mathbb{R}, \ l(g) = a, \ l(g^2) = b,$$

$$\text{where } f_t(s) = I(s \leq t), g(s) = s, s \in \mathbb{R}\Big\}$$

**Theorem 2.3.** *Let $\psi : \mathbb{R} \to [0,\infty)$ be a function such that $\int_\mathbb{R} \psi(x)\,dx < \infty$, let*

$$\text{(2.8)} \quad H_n^{\text{AD}}(u) := P_\Phi\left\{\int_{-\infty}^{\infty}(F_n(\bar{X}_n + s_n t) - \Phi(t))^2\psi(\Phi(t))d\Phi(t) \geq u\right\}, u \geq 0,$$

*and let*

$$G^{\text{AD}}(u) := \inf\Big\{J_\Phi(l) : l \in (\mathcal{L}_\Phi^{\Upsilon_0})^*,$$

(2.9)
$$\int_{-\infty}^{\infty}[x(a + (b-a^2)^{1/2}t) - \Phi(t)]^2\psi(F_0(t))dF_0(t) \geq u,$$

$$l(f_t) = x(t), t \in \mathbb{R}, \ l(g) = a, \ l(g^2) = b,$$

$$\text{where } f_t(s) = I(s \leq t), g(s) = s, s \in \mathbb{R}\Big\}.$$

*Let $\{X_j\}_{j=1}^{\infty}$ be a sequence of i.i.d.r.v.'s with a continuous c.d.f. $F$. Then,*

$$-\lim_{\delta\to 0+} G^{\text{AD}}\left(\int_{-\infty}^{\infty}[F(\mu_F + \sigma_F t) - \Phi(t)]^2\psi(\Phi(t))d\Phi(t) + \delta\right)$$

(2.10)
$$\leq \liminf_{n\to\infty} n^{-1}\ln H_n^{\text{AD}}\left(\int_{-\infty}^{\infty}[F_n(\bar{X}_n + s_n t)) - \Phi(t)]^2\psi(\Phi(t))d\Phi(t) > u\right)$$

$$\leq \limsup_{n\to\infty} n^{-1}\ln H_n^{\text{AD}}\left(\int_{-\infty}^{\infty}[F_n(\bar{X}_n + s_n t) - \Phi(t)]^2\psi(\Phi(t))d\Phi(t) \geq u\right)$$

$$\leq -\lim_{\delta\to 0+} G^{\text{AD}}\left(\int_{-\infty}^{\infty}[F(\mu_F + \sigma_F t) - \Phi(t)]^2\psi(\Phi(t))d\Phi(t) - \delta\right) \text{ a.s.}$$

## 3. Simulations

We present simulations of the mean of the p-value of several alternatives. As before, suppose that a test rejects $H_0$ if $T_n \geq c$, where $T_n := T_n(X_1, \ldots, X_n)$ is a statistic and $c$ is a constant. The significance level of the test is

$$(3.1) \qquad\qquad H_n(c) := \sup_{\theta \in \Theta_0} P_\theta(T_n \geq c),$$

where $P_\theta$ denotes the probability measure for which the data has p.d.f. $f(\cdot, \theta)$. The $p$–value of the test is $H_n(T_n)$. We do simulations estimating $E[H_n(T_n)]$.

Let $T_n^1, \ldots, T_n^N$ be $N$ simulations of the test under the null hypothesis using a sample size $n$. Let $T_n^{1,\text{alt}}, \ldots, T_n^{1,\text{alt}}$ be $N$ simulations of the test under a certain alternative hypothesis. Then, $N^{-2} \sum_{j,k=1}^N I(T_n^j \geq T_n^{k,\text{alt}})$ estimates $E[H_n(T_n)]$, where the expectation is taken assuming that $T_n$ is obtained using $n$ i.i.d.r.v.s from the alternative hypothesis. In Table 1, $N = 10000$ is used for the Lilliefors, the Cramer–von Mises, the Anderson–Darling, the Shapiro-Wilk, and the BHEP test ([11] and [7]).

TABLE 1

| $n$ | L | CM | AD | SW | BHEP |
|---|---|---|---|---|---|
| Alternative: exponential distribution | | | | | |
| 10 | 0.248325 | 0.2065738 | 0.1878327 | 0.1621557 | 0.1813481 |
| 15 | 0.1601991 | 0.1178508 | 0.0946044 | 0.07611985 | 0.09569895 |
| 20 | 0.1043946 | 0.06510648 | 0.05291726 | 0.03304067 | 0.05206452 |
| 30 | 0.044566 | 0.02152872 | 0.0129459 | 0.00750638 | 0.01409681 |
| 50 | 0.00818707 | 0.00203949 | 0.0009082 | 0.00241882 | 0.00121646 |
| Alternative: double exponential distribution | | | | | |
| 10 | 0.3992631 | 0.3939724 | 0.3983314 | 0.4148731 | 0.4109459 |
| 15 | 0.3499758 | 0.339891 | 0.3389608 | 0.3677549 | 0.3640368 |
| 20 | 0.3169975 | 0.2979569 | 0.3009778 | 0.3294354 | 0.3244178 |
| 30 | 0.2616672 | 0.2341172 | 0.2397223 | 0.2807123 | 0.2555247 |
| 50 | 0.1796492 | 0.1417312 | 0.1434135 | 0.2837385 | 0.1564005 |
| Alternative: Cauchy distribution | | | | | |
| 10 | 0.1566254 | 0.1493284 | 0.1500601 | 0.1705618 | 0.1682588 |
| 15 | 0.08474307 | 0.07479607 | 0.07505173 | 0.09128725 | 0.08964657 |
| 20 | 0.04651569 | 0.03862244 | 0.03767999 | 0.04857044 | 0.04194876 |
| 30 | 0.01420118 | 0.0100633 | 0.00974881 | 0.01496182 | 0.01179398 |
| 50 | 0.0017361 | 0.00066862 | 0.00087474 | 0.00405 | 0.00048095 |
| Alternative: Beta(2,1) distribution | | | | | |
| 10 | 0.41099 | 0.3884071 | 0.3686801 | 0.3358273 | 0.3565215 |
| 15 | 0.3608343 | 0.321224 | 0.2986805 | 0.2631936 | 0.2861669 |
| 20 | 0.3147082 | 0.272089 | 0.2446055 | 0.1861695 | 0.2330953 |
| 30 | 0.2411935 | 0.1840438 | 0.1534217 | 0.09638084 | 0.1428481 |
| 50 | 0.1355312 | 0.08707787 | 0.05502258 | 0.0835284 | 0.0572854 |
| Alternative: Beta(3,3) distribution | | | | | |
| 10 | 0.5084849 | 0.5040435 | 0.51387 | 0.4928259 | 0.4947155 |
| 15 | 0.5063525 | 0.5029377 | 0.5033103 | 0.4872658 | 0.4864098 |
| 20 | 0.5072011 | 0.5037995 | 0.4993457 | 0.4762998 | 0.4797991 |
| 30 | 0.4899722 | 0.4745532 | 0.4780857 | 0.4285843 | 0.4510982 |
| 50 | 0.4590308 | 0.4447785 | 0.4382911 | 0.4189339 | 0.4064414 |

Table 1 (Continued)

Alternative: Logistic(1) distribution

| | | | | | |
|----|-----------|-----------|-----------|-----------|-----------|
| 10 | 0.4736219 | 0.4725762 | 0.4685902 | 0.4749748 | 0.4677617 |
| 15 | 0.4560905 | 0.4468502 | 0.4335687 | 0.4624648 | 0.46532 |
| 20 | 0.4493409 | 0.4488339 | 0.4450634 | 0.4426041 | 0.4510982 |
| 30 | 0.4410423 | 0.422886 | 0.4233825 | 0.4348024 | 0.4153006 |
| 50 | 0.4204326 | 0.3978938 | 0.3770672 | 0.4458524 | 0.3819914 |

Alternative: uniform distribution

| | | | | | |
|----|-----------|-----------|-----------|-----------|------------|
| 10 | 0.4438842 | 0.4241476 | 0.404153 | 0.3681328 | 0.3922898 |
| 15 | 0.4059967 | 0.3716177 | 0.3468994 | 0.3023148 | 0.338187 |
| 20 | 0.3739766 | 0.3308353 | 0.2951247 | 0.2270906 | 0.27552 |
| 30 | 0.3050368 | 0.2429758 | 0.2024329 | 0.117917 | 0.1862224 |
| 50 | 0.2066687 | 0.1359771 | 0.0889871 | 0.1007488 | 0.08932786 |

We should expect the average $p$–value is small than 0.5. However, for the Beta(3,3) distribution the average $p$–value is bigger than 0.5. Between the tree test which use the empirical c.d.f., the Anderson–Darling test is the one with smallest average $p$–value. For almost half of the considered distributions, this is the test with the smallest average $p$–value. The Shapiro–Wilk test also performs very well overall.

## 4. Proofs

We will need the following lemmas:

**Lemma 4.1** (Lemma 5.1 (i), [4]). *For each $k \geq 0$ and each function $f \in \mathcal{L}^{\Upsilon_0}$,*

$$|l(f)| \leq (J(l) + 1 + 2^{1/2}) N_{\Upsilon_0}(f).$$

**Lemma 4.2.** *Let $l \in (\mathcal{L}_{\Phi}^{\Upsilon_0})^*$ with $J_{\Phi}(l) < \infty$. Then, $x(t) = l(I(\cdot \leq t))$, $t \in \mathbb{R}$, is a continuous function with $\lim_{t \to -\infty} x(t) = 0$ and $\lim_{t \to \infty} x(t) = 1$.*

*Proof.* Let $\alpha : (0, \infty) \to (0, \infty)$ be defined as $\alpha(x) = \exp(1/x) - 1/x$. It is easy to see that $\alpha$ is one–to–one function. We claim that for a Borel set $A \subset \mathbb{R}$,

$$(4.1) \qquad N_{\Upsilon_0}(I(X \in A)) = \alpha^{-1}(1 + (P_{\Phi}(X \in A))^{-1}),$$

where $\alpha^{-1}$ denotes the inverse function of $\alpha$. We have that

$$E_{\Phi}[\Upsilon_0(\lambda^{-1} I(X \in A))] = E[\exp(\lambda^{-1} I(X \in A)) - 1 - \lambda^{-1} I(X \in A)]$$
$$= p \exp(\lambda^{-1}) + 1 - p - 1 - \lambda^{-1} p,$$

where $p := P_{\Phi}(X \in A)$. So, $1 \geq E_{\Phi}[\Upsilon_0(\lambda^{-1} I(X \in A))]$ is equivalent to $\alpha(\lambda) \leq 1 + p^{-1}$. So, (4.1) follows.

By Lemma 4.1 and (4.1), for each $s, t \in \mathbb{R}$ with $s < t$,

$$|x(t)| = |l(I(X \leq t))| \leq (J_{\Phi}(l) + 1 + 2^{1/2}) \alpha^{-1}(1 + (\Phi(t))^{-1}),$$

$$|1 - x(t)| = |l((X > t))| \leq (J_{\Phi}(l) + 1 + 2^{1/2}) \alpha^{-1}(1 + (1 - \Phi(t))^{-1}),$$

and

$$|x(t) - x(s)| \leq (J_{\Phi}(l) + 1 + 2^{1/2}) \alpha^{-1}(1 + (\Phi(t) - \Phi(s))^{-1}),$$

which implies the claim. □

*Proof of Lemma 2.1.* By Theorem 2.1, $\{U_n(t) : t \in \mathbb{R}\}$ satisfies the LDP in $l_\infty(\mathbb{R})$ with speed $n$, where $U_n(t) = F_n(t)$. Let $\omega_1$ and $\omega_2$ be two numbers, which are not in $\mathbb{R}$. Let $U_n(\omega_1) = n^{-1} \sum_{j=1}^n X_j$ and let $U_n(\omega_2) = n^{-1} \sum_{j=1}^n X_j^2$. By the LDP for sums of i.i.d. $\mathbb{R}^d$–valued r.v.'s, the finite dimensional distributions of $\{U_n(t) : t \in \mathbb{R} \cup \{\omega_1, \omega_2\}\}$ satisfy the LDP with speed $n$. Since $\{U_n(t) : t \in \mathbb{R}\}$ satisfies the LDP in $l_\infty(\mathbb{R})$, it satisfies an exponential asymptotic equicontinuity condition (see Theorem 3.1 in [2]). This implies that $\{U_n(t) : t \in \mathbb{R} \cup \{\omega_1, \omega_2\}\}$ satisfies an exponential asymptotic equicontinuity condition. So, $\{U_n(t) : t \in \mathbb{R} \cup \{\omega_1, \omega_2\}\}$ satisfies the LDP in $l_\infty(\mathbb{R} \cup \{\omega_1, \omega_2\})$ with speed $n$ (see Theorem 3.1 in [2]). Besides, the rate of function is

$$I(x) = \inf\{J_\Phi(l) : l \in (\mathcal{L}^{\Upsilon_0})^*, l(I(X \leq t) = x(t), t \in \mathbb{R}, l(X) = x(\omega_1),$$
$$l(X^2) = x(\omega_2), x \in l_\infty(\mathbb{R} \cup \{\omega_1, \omega_2\})\}$$

Let $\Gamma_n : l_\infty(\mathbb{R} \cup \{\omega_1, \omega_2\}) \to \mathbb{R}$ be defined by

$$\Gamma_n(x) = \sup_{t \in \mathbb{R}} |x(x(\omega_1) + n^{1/2}(n-1)^{-1/2}(x(\omega_2) - (x(\omega_1))^2)^{1/2}t) - \Phi(t)|\psi(t),$$

for $x \in l_\infty(\mathbb{R} \cup \{\omega_1, \omega_2\})$. Next, we prove using Theorem 2.1 in Arcones (2003a) that

$$\Gamma_n(\{U_n(t) : t \in \mathbb{R} \cup \{\omega_1, \omega_2\}\}) = \sup_{t \in \mathbb{R}} |F_n(\bar{X}_n + s_n t) - \Phi(t)|\psi(t)$$

satisfies the LDP in $\mathbb{R}$ with speed $n$ and rate function

$$(4.2) \quad \begin{aligned} q^{\mathrm{Li}}(u) :&= \inf\{J_\Phi(l) : l \in (\mathcal{L}^{\Upsilon_0})^*, \sup_{t \in \mathbb{R}} |x(a + (b - a^2)^{1/2}t) - \Phi(t)|\psi(t) = u, \\ &\quad l(I(X \leq t) = x(t), t \in \mathbb{R}, \ l(X) = a, \ l(X^2) = b\}. \end{aligned}$$

To apply Theorem 2.1 in [2], we need to prove that if $x_n \to x$, in $l_\infty(\mathbb{R} \cup \{\omega_1, \omega_2\})$ and $I(x) < \infty$, then $\Gamma_n(x_n) \to \Gamma(x)$, where

$$\Gamma(x) = \sup_{t \in \mathbb{R}} |x(x(\omega_1) + (x(\omega_2) - (x(\omega_1))^2)^{1/2}t) - \Phi(t)|\psi(t), x \in l_\infty(\mathbb{R} \cup \{\omega_1, \omega_2\}).$$

We have that

$$\begin{aligned} |\Gamma_n(x_n) &- \Gamma(x)| \\ &\leq \sup_{t \in \mathbb{R}} |x_n(x_n(\omega_1) + n^{1/2}(n-1)^{-1/2}(x_n(\omega_2) - (x_n(\omega_1))^2)^{1/2}t)\psi(t) \\ &\quad - x(x_n(\omega_1) + n^{1/2}(n-1)^{-1/2}(x_n(\omega_2) - (x_n(\omega_1))^2)^{1/2}t)\psi(t)| \\ &\quad + \sup_{t \in \mathbb{R}} |x(x_n(\omega_1) + n^{1/2}(n-1)^{-1/2}(x_n(\omega_2) - (x_n(\omega_1))^2)^{1/2}t)\psi(t) \\ &\quad - x(x(\omega_1) + (x(\omega_2) - (x(\omega_1))^2)^{1/2}t)\psi(t)| \\ &=: I + II. \end{aligned}$$

Since $x_n \to x$, in $l_\infty(\mathbb{R} \cup \{\omega_1, \omega_2\})$,

$$I = \sup_{t \in \mathbb{R}} |x_n(t) - x(t)|\psi(t) \to 0.$$

By Lemma 4.2, $x$ is a continuous function with $\lim_{t \to -\infty} x(t) = 0$ and $\lim_{t \to \infty} x(t) = 1$. So,

$$II \to 0.$$

From the previous computations, we get that $\Gamma_n(x_n) \to \Gamma(x)$. Hence, $\sup_{t \in \mathbb{R}} |F_n \times (\bar{X}_n + s_n t) - \Phi(t)|$ satisfies the LDP in $\mathbb{R}$ with speed $n$ and rate function $q^{\mathrm{Li}}(u)$. This implies (2.4).  □

*Proof of Theorem 2.2.* We have that

$$|\sup_{t\in\mathbb{R}}|F_n(\bar{X}_n + s_n t)) - \Phi(t)|\psi(t) - \sup_{t\in\mathbb{R}}|F(\mu_F + \sigma_F t) - \Phi(t)|\psi(t)|$$

$$\leq \sup_{t\in\mathbb{R}}|F_n(\bar{X}_n + s_n t) - F(\mu_F + \sigma_F t)|$$

$$\leq \sup_{t\in\mathbb{R}}|F_n(t) - F(t)|\psi(t) + \sup_{t\in\mathbb{R}}|F(\bar{X}_n + s_n t)) - F(\mu_F + \sigma_F t)|\psi(t).$$

By the Glivenko–Cantelli theorem,

$$\sup_{t\in R}|F_n(t) - F(t)|\psi(t) \to 0 \quad a.s.$$

Using that

$$\bar{X}_n \to \mu_F \quad a.s., s_n \to \sigma_F \quad a.s.$$

and $F$ is a continuous function with $\lim_{t\to-\infty}F(t) = 0$ and $\lim_{t\to\infty}F(t) = 1$, we get that

$$\sup_{t\in\mathbb{R}}|F(\bar{X}_n + s_n t) - F(\mu_F + \sigma_F t)|\psi(t) \to 0 \quad a.s.$$

Hence,

$$(4.3) \qquad \sup_{t\in\mathbb{R}}|F_n(\bar{X}_n + s_n t) - \Phi(t)|\psi(t) \to \sup_{t\in\mathbb{R}}|F(\mu_F + \sigma_F t) - \Phi(t)|\psi(t)| \quad a.s.$$

The claim in this theorem follows from (4.3) and Lemma 2.1. □

The proofs of Lemma 2.2 and Theorem 2.3 are similar to those of Lemma 2.1 and Theorem 2.2 and they are omitted.